\documentclass[a4paper,10pt]{article}
\usepackage{amsthm}
\usepackage{amsmath}
\usepackage{amsfonts}
\usepackage{amssymb}
\usepackage{latexsym}
\usepackage[margin=1.2in]{geometry}

\usepackage{graphicx}

\usepackage[table]{xcolor}

\usepackage{tabularx}

\newtheorem{prop}{Proposition}[section]
\newtheorem{thm}[prop]{Theorem}
\newtheorem{cor}[prop]{Corollary}

\newtheorem{rem}[prop]{Remark}

\newtheorem*{KN}{Crazy Knight's Tour Problem}
\newtheorem{lem}[prop]{Lemma}

\newtheorem{ex}[prop]{Example}

\newtheorem{defin}[prop]{Definition}

\numberwithin{equation}{section}

\newcommand{\probname}{Crazy Knight's Tour Problem}

\title{A tour problem on a toroidal board}
\date{}

\author{Simone Costa
\thanks{DII/DICATAM - Sez. Matematica, Universit\`a degli Studi di Brescia, Via
Branze 38, I-25123 Brescia, Italy. email: simone.costa@unibs.it}\ \
Marco Dalai
\thanks{DII, Universit\`a degli Studi di Brescia, Via
Branze 38, I-25123 Brescia, Italy. email: marco.dalai@unibs.it} \ \
 and Anita Pasotti
 \thanks{DICATAM - Sez. Matematica, Universit\`a degli Studi di Brescia, Via
Branze 43, I-25123 Brescia, Italy. email: anita.pasotti@unibs.it}}

\begin{document}
\maketitle

\begin{abstract}
In this paper we study a tour problem that we came cross while studying biembeddings and Heffter arrays,
see [D.S. Archdeacon,
Heffter arrays and biembedding graphs on surfaces,
\textit{Electron. J. Combin.} 22 (2015) \#P1.74].
Let $A$ be an $n\times m$ toroidal array consisting of filled cells and empty cells. Assume that an orientation $R=(r_1,\dots,r_n)$
of each row and $C=(c_1,\dots,c_m)$ of each column of $A$ is fixed. Given an initial filled cell $(i_1,j_1)$ consider the list
$ L_{R,C}=((i_1,j_1),(i_2,j_2),\ldots,(i_k,j_k),$ $(i_{k+1},j_{k+1}),\ldots)$
where $j_{k+1}$ is the column index of the filled cell $(i_k,j_{k+1})$ of the row $R_{i_k}$ next to $(i_k,j_k)$ in the orientation $r_{i_k}$,
and where $i_{k+1}$ is the row index of the filled cell of the column $C_{j_{k+1}}$ next to $(i_k,j_{k+1})$ in the orientation $c_{j_{k+1}}$.
 We propose the following
\begin{KN}
Do there exist $R$ and $C$ such that the list $L_{R,C}$ covers all the filled cells of $A$?
\end{KN}
Here we provide a complete solution for the case with no empty cells and we obtain partial results for square arrays where the
filled cells follow some specific regular patterns.
\end{abstract}

\section{Introduction}
Chessboards have been a classical setting for many challenging problems starting from very ancient times up to the most recent ones, for a survey see \cite{HHR,W}.
Most of these problems are not only interesting \emph{per se}, but have many applications
and often arise from a mathematical context.
For example classical problems of covering chessboards with the minimum number of chess pieces are related
to the study of dominating sets in graphs, see for instance \cite{C,CH}.
As well as domination problems, even tour problems are much studied.
 The most famous problem of this type is the \emph{Knight's tour problem} (see for instance \cite{CHMW}),
whose goal is to determine a series of moves made by a knight so that it visits every square on a chessboard exactly once.
Then several variations of the original statement have been proposed and investigated, see \cite{BYZJH,CO,EGG,K,MF,MF2} and the references therein.

Our aim is to find out how to obtain a closed tour on an $n \times m$ board
wrapped onto a torus so that the $n$ rows and the $m$ columns go around the torus.
Hence, in the following we will consider all the row indices modulo $n$ and all the column indices modulo $m$.
Also, the general board we consider might have some squares removed, and we thus represent it as a matrix $A$ with some empty
and some filled cells, naming it a
 \emph{partially filled array}. The tour has to visit every filled cell exactly once.
Since we are only interested in which cells are empty and which filled, but not in the elements of $A$, we will denote the filled positions simply with a $\bullet$.

In order to state our problem we  introduce a \emph{move function} which defines the rules for the considered tour on the board.
Given two positive integers $a$ and $b$ with $a\leq b$, by $[a,b]$ we will mean the set $\{a,a+1,\ldots,b\}$.
Let us consider a partially filled $n\times m$ array $A$; we denote by $F(A)$ the subset of $[1,n]\times [1,m]$ given by the filled positions of $A$.
Given  $(i,j)\in F(A)$, we define the row successor $s_r((i,j))$ to be $(i,j+k)$ where $k\geq 1$ is the minimum value such that $(i,j+k)\in F(A)$.
Similarly we define the column successor $s_c((i,j))$ to be $(i+k,j)$ where $k\geq1$ is the minimum value such that $(i+k,j)\in F(A)$.
Given two vectors $R:=(r_1,\dots,r_n)\in \{-1,1\}^n$ and $C:=(c_1,\dots,c_m)\in \{-1,1\}^m$, we call \emph{move function} the function $S_{\tiny{R,C}}: F(A)\rightarrow F(A)$ so defined
$$
S_{R,C}((i,j)):=
\begin{cases}
 s_c(s_r((i,j))) \mbox{ if } r_i=1 \mbox{ and } c_{j'}=1 \mbox{ where } s_r((i,j))= (i,j');\\
 s_c(s_r^{-1}((i,j))) \mbox{ if } r_i=-1 \mbox{ and } c_{j'}=1 \mbox{ where } s_r^{-1}((i,j))= (i,j');\\
 s_c^{-1}(s_r((i,j))) \mbox{ if } r_i=1 \mbox{ and } c_{j'}=-1 \mbox{ where } s_r((i,j))= (i,j');\\
 s_c^{-1}(s_r^{-1}((i,j))) \mbox{ if } r_i=-1 \mbox{ and } c_{j'}=-1 \mbox{ where } s_r^{-1}((i,j))= (i,j').\\
\end{cases}
$$
We can see the vector $R$ as a choice of the direction of each row: from left to right if $r_i=1$ and from right to left if $r_i=-1$.
Similarly the vector $C$ can be seen as a choice of the direction of each column: from top to bottom if $c_i=1$ and in the reverse way if $c_i=-1$.
Moreover, we can give the following interpretation to the move $S_{R,C}((i,j))$: from the position $(i,j)$ we move first in the $i$-th row following the direction of $r_i$ and then, from the arrival position $(i,j')$ we move in the $j'$-th column in direction of $c_{j'}$ arriving in the position $S_{R,C}((i,j))$.
For any $(i,j)\in F(A)$, we set $$L(i,j):=((i,j),S_{R,C}((i,j)), S_{R,C}^2((i,j)),\dots, S_{R,C}^{p}((i,j))),$$
where $p$ is the minimum positive integer such that $S_{R,C}^{p+1}((i,j))=(i,j)$,
namely $L(i,j)$ is the tour that we obtain starting from the cell $(i,j)$. It is natural to ask whether, in this way, we can cover all the filled positions or not.

\begin{KN}\label{Arco}
Given a partially filled $n\times m$ array $A$, determine whether there exist vectors $R\in \{-1,1\}^n$ and $C\in \{-1,1\}^m$ such that, given a filled position $(i,j)$, the list
$L(i,j)$ covers all the filled positions of $A$.
\end{KN}
By $P(A)$ we will denote the \probname\ for a given array $A$.
Also, given a filled cell $(i,j)$, if $L(i,j)$ covers all the filled positions of $A$ we will say that the vectors $R$ and $C$ are a solution of $P(A)$.

In this paper, firstly,  we will present some necessary conditions for the existence of a solution of $P(A)$ where $A$ is a given array, see Section \ref{sec2}.
Then, in Section \ref{sec:piene}, we present a complete solution when $A$ is a totally filled rectangular array.
We focus also on square arrays with some empty cells.
In particular in Section \ref{SecDiag} and in Section \ref{SecAlmost} we obtain partial results when $A$ has exactly $k$ filled cells in each row and in each column and
when $A$ has exactly $k$ filled cells in each row and in each column except for one row and column which have $k+1$ filled positions, respectively.
Finally, in Section \ref{SecRecursive}, we present a recursive construction which allows us to obtain other infinite classes of arrays $A$
such that $P(A)$ has a solution.

\subsection{Motivation}
In \cite{A}, Archdeacon introduced the concept of a Heffter array as a useful tool for determining biembeddings
of complete graphs, that is $2$-colorable embeddings. In details, in that paper the author investigated the case in which
the face boundaries of the two colour classes form two cycle systems, for a survey see \cite{GG}.
Heffter arrays which give rise to biembeddings have been constructed in \cite{CMPPGlobally, CMPPRelative, DM}.
In particular, in \cite{CMPPGlobally} the authors introduced the class of \emph{globally simple} Heffter arrays
and, implicitly, studied the \probname\ in some special instances (see \cite[Propositions 3.4 and 3.6]{CMPPGlobally})
obtaining in such a way new biembeddings (see \cite[Theorem 1.11]{CMPPGlobally}).
The relationship between globally simple Heffter arrays, \probname\ and biembeddings is explained in the following result,
that is a reformulation of  \cite[Theorem 1.1]{A}, in the case of square globally simple Heffter arrays, in terms of $P(A)$. Clearly, a similar theorem holds in the rectangular case.
\begin{thm}
  Let $A$ be a globally simple Heffter array $n\times n$ such that each row and each column has exactly $k$ filled cells.
   If there exists a solution of $P(A)$, then there exists a biembedding of the complete graph of order $2nk+1$
   on orientable surface whose face boundaries are $k$-cycles.
\end{thm}
We point out that Heffter arrays are considered interesting as combinatorial objects on their own
and not only in view of their relationship with biembeddings. For instance, in \cite{ABD, ADDY, BCDY, CDDY,RelH, DW}
the authors construct new infinite classes of Heffter arrays without investigating biembeddings.
Analogously, we believe that also the \probname\ is interesting \emph{per se}, as many other tour problems
which have been largely studied.

\section{Preliminary considerations}\label{sec2}
Clearly, we can suppose that each row and each column of $A$ has at least one filled position.
We note that if the $i$-th row $R_i$ or the $j$-th column $C_j$ has exactly just one filled position $(i,j)$
then $s_r((i,j))=(i,j)$ or $s_c((i,j))=(i,j)$, respectively.
Also if there exist a row $R_i$ and a column $C_j$ such that $R_i\cup C_j$ has exactly one filled position,
then this position is $(i,j)$, since we have assumed that we have no empty row or empty column.
In this case, since $S_{R,C}((i,j))=(i,j)$, it is immediate that the \probname\ has no solution, except for the trivial case in which $A$ is a square array of size $1$.
\begin{rem}\label{11=ij}
Given an array $A$ and a filled position $(i,j)$,
the list $L(i,j)$ covers all the filled positions of $A$ if and only if the list $L(i',j')$ covers all the filled positions of $A$, for any $(i',j')\in F(A)$. Since we are dealing with toroidal boards, for the aim of this paper it is not restrictive to suppose that $(1,1)\in F(A)$, namely that the cell $(1,1)$ of the given array $A$ is filled.

\end{rem}

\begin{ex}\label{ex0}
Let $A$ be
$$\begin{array}{|r|r|r|r|}\hline
     \bullet &  & & \bullet      \\\hline
      &   \bullet &   \bullet &      \\\hline
      &   \bullet  &  \bullet & \bullet    \\\hline
      \bullet &    &    \bullet  &    \\\hline
    \end{array}$$
		Choosing $R:=(-1,1,1,-1)$ and $C:=(1,-1,1,1)$ we can cover all the filled positions of $A$, as shown in the table below where
		in each filled position we write $j$ if we reach that position after having applied $S_{R,C}$ to $(1,1)$ exactly $j$ times.
Here we represent the elements of $R$ and $C$ by an arrow.		
$$\begin{array}{r|r|r|r|r|}
  &  \downarrow   & \uparrow & \downarrow & \downarrow   \\\hline
\leftarrow &     0 &      &    &  4  \\\hline
\rightarrow &      &   2 &   6 &      \\\hline
\rightarrow &      &   7  & 3  & 1   \\\hline
\leftarrow &      5 &    &  8  &    \\\hline
    \end{array}$$
\end{ex}

\begin{rem}
Given a partially filled $n\times m$ array $A$, let $R:=(r_1,r_2,\ldots,r_n)$ and $C:=(c_1,c_2,\ldots,\\c_m)$ be
 a solution of $P(A)$. Note that this does not imply that $R^i:=(r_i,\ldots,r_n,r_1,\ldots,r_{i-1})$ and $C^j:=(c_j,\ldots,c_m,c_1,\ldots,c_{j-1})$
are a solution of $P(A)$ too, where $i=2,\ldots,n$ and $j=2,\ldots,m$.
\end{rem}
\begin{ex}
Let $A$ be the array of Example \ref{ex0}. If instead of $R=(-1,1,1,-1)$ and $C=(1,-1,1,1)$ we take $R^2=(1,1,-1,-1)$ and $C^2=(-1,1,1,1)$
it is easy to see that we do not cover all the filled cells of $A$.
\end{ex}

In order to determine the necessary conditions for the existence of a solution of the \probname, we introduce the notion of
\emph{closed subarray}.
\begin{defin} Let $A$ be an $n \times m$ toroidal array with no empty row or column. Let $\mathcal{R}$ be a list of rows of $A$ and let $\mathcal{C}$ be a list of columns of $A$.
We say that the subarray $\mathcal{R}\cap \mathcal{C}$ of $A$ is \emph{closed}  if  $F(\mathcal{R}\cap\mathcal{C})=F(\mathcal{R}\cup \mathcal{C})$.
We say that a closed subarray is minimal if it is
minimal with respect to the inclusion.
\end{defin}

\begin{ex}
Let $A$ be the following $7\times 9$ array.
$$\begin{array}{|r|r|r|r|r|r|r|r|r|}\hline
 &   &   &    \bullet &  &     &    &    &   \\ \hline
  \bullet &    &  \bullet &  & \bullet      &   &  \bullet    &  \bullet    & \bullet      \\ \hline
 &   &   &  \bullet &    &   \bullet  &    &    &     \\ \hline
  \bullet &   \bullet  &   \bullet  & &  \bullet   &   &  \bullet    &  \bullet    & \bullet      \\ \hline
  \bullet &   \bullet  &   \bullet  & &  \bullet   &   &  \bullet    &      &      \\ \hline
 &   &   &  \bullet    &&   \bullet  &    &    &    \\ \hline
 &   &   &    \bullet & &  \bullet   &    &    &    \\ \hline
    \end{array}$$
Consider the lists $\mathcal{R}=(R_2,R_4,R_5)$  and  $\mathcal{C}=(C_1,C_2,C_3,C_5,C_{7},C_8,C_{9})$.
Then it is easy to see that $\mathcal{R}\cap \mathcal{C}$ is a closed subarray of $A$.
\end{ex}

\begin{thm}\label{prop:necc}
Let $A$ be an $n\times m$ array with $|F(A)|$ filled cells. Necessary conditions for the existence of a solution of $P(A)$ are:
\begin{itemize}
  \item[1)] the array $A$ is a minimal closed subarray;
  \item[2)] $|F(A)|\equiv m+n-1 \pmod 2$.
\end{itemize}
\end{thm}
\proof
1) Let $\mathcal{R}\cap\mathcal{C}$ be a closed subarray of $A$ such that $\mathcal{R}\cap \mathcal{C}\not=A$. Since $A$ has no empty row and no empty column, $F(\mathcal{R}\cap \mathcal{C}) \subset F(A)$.
It is immediate to see that for any $(i,j)\in F(\mathcal{R}\cap \mathcal{C})$, $L(i,j)\subseteq F(\mathcal{R}\cap \mathcal{C})$, hence it does not cover all the filled positions of $A$.

2) The function $S_{R,C}$ defines a permutation $\omega$ on $F(A)$.
Let $\omega_r$ and $\omega_c$ be the permutations on the rows and on the columns, respectively, obtained from the move function.
By definition of $S_{R,C}$, we have $\omega=\omega_c \circ \omega_r$.
Hence the parity of $\omega$, that is $|F(A)|-1$, has to be equal to the parity of $\omega_c \circ \omega_r$, that is
$\sum_{i=1}^{n}(|F(R_i)|-1)+\sum_{j=1}^{m}(|F(C_j)|-1)=\sum_{i=1}^{n}|F(R_i)|+\sum_{j=1}^{m}|F(C_j)|-n-m=2|F(A)|-n-m$.
It follows that  $|F(A)|\equiv n+m-1 \pmod 2$.
\endproof

\begin{cor}\label{cor:necc}
Let $A$ be a square array with $|F(A)|$ filled cells. If there exists a solution of $P(A)$ then $|F(A)|$ is odd.
\end{cor}

\section{Totally filled arrays}\label{sec:piene}
In this section, we solve Problem $P(A)$ when $A$ is a totally filled $n \times m$ array (i.e. $|F(A)|=mn$).
In this case the move function is simply given by:
$$
S_{R,C}((i,j)):=
\begin{cases}
 (i+1,j+1) \mbox{ if } r_i=1 \mbox{ and } c_{j+1}=1;\\
 (i+1,j-1) \mbox{ if } r_i=-1 \mbox{ and } c_{j-1}=1;\\
 (i-1,j+1) \mbox{ if } r_i=1 \mbox{ and } c_{j+1}=-1;\\
 (i-1,j-1) \mbox{ if } r_i=-1 \mbox{ and } c_{j-1}=-1.
\end{cases}
$$
We first introduce a useful lemma.
\begin{lem}\label{Coprimi}
Let $A$ be a totally filled $n \times m$ array.
Let $R:=(1,1,\dots,1)$ and $C:=(1,\dots,c_{m-l}=1,c_{m-l+1}=-1,\dots,-1)$, where $0\leq l\leq m$. Then $R$ and $C$ are a solution of $P(A)$ if and only if $\gcd(m-2l,n)=1$.
\end{lem}
\proof
Note that, for $i=1,\ldots,n$, the first $m$ cells of $L(i,1)$ form the list
\begin{multline}
 \Gamma_i:=((i,1),(i+1,2),\dots, (i+m-l-1, m-l),\\
  (i+m-l-2, m-l+1),(i+m-l-3, m-l+2),\ldots, (i+m-2l-1,m))\,.
 \end{multline}
Since the $\Gamma_i$'s contain one element in each column and are vertical shift of one another, it is clear that
$$\bigcup_{i=1}^{n}\Gamma_i=F(A)\,.
$$
Now, since $c_1=1$, $L(1,1)=(\Gamma_1, \Gamma_{1+(m-2l)}, \Gamma_{1+2(m-2l)},\ldots)$,
where the indices are taken modulo $n$. The indices cover all the values $1, 2, \ldots, n$ if and only if $\gcd(m-2l,n)=1$, so this is equivalent to the condition $L(1,1)=F(A)$.
\endproof
\begin{ex}
Consider the totally filled array $A$ of size $n\times m$ with $n=5$ and $m=14$. Also let $l=3$. In the table below, cells in  $\Gamma_i$ are flagged with $\gamma_i$.
\begin{center}
\begin{footnotesize}
$\begin{array}{|r|r|r|r|r|r|r|r|r|r|r|r|r|r|}\hline
\gamma_1 &  \gamma_5 &  \gamma_4 &  \gamma_3  & \gamma_2  & \gamma_1   & \gamma_5   &  \gamma_4 &  \gamma_3 & \gamma_2 & \gamma_1 & \gamma_2 & \gamma_3 & \gamma_4\\ \hline
\gamma_2 & \gamma_1 &  \gamma_5 &  \gamma_4 &  \gamma_3  & \gamma_2  & \gamma_1   & \gamma_5   &  \gamma_4 &  \gamma_3 & \gamma_2  & \gamma_3 & \gamma_4 & \gamma_5\\ \hline
\gamma_3 & \gamma_2 & \gamma_1 &  \gamma_5 &  \gamma_4 &  \gamma_3  & \gamma_2  & \gamma_1   & \gamma_5   &  \gamma_4 &  \gamma_3  & \gamma_4 & \gamma_5 & \gamma_1\\ \hline
\gamma_4 & \gamma_3 & \gamma_2 & \gamma_1 &  \gamma_5 &  \gamma_4 &  \gamma_3  & \gamma_2  & \gamma_1   & \gamma_5   &  \gamma_4  & \gamma_5 & \gamma_1 & \gamma_2\\ \hline
\gamma_5 & \gamma_4 & \gamma_3 & \gamma_2 & \gamma_1 &  \gamma_5 &  \gamma_4 &  \gamma_3  & \gamma_2  & \gamma_1   & \gamma_5   & \gamma_1 & \gamma_2 & \gamma_3\\ \hline
    \end{array}$
		\end{footnotesize}
				\end{center}
				Since $\gcd(m-2l,n)=\gcd(8,5)=1$ we can apply Lemma \ref{Coprimi}.
So we start from the position $(1,1)$ and apply $S_{R,C}$ until we arrive again in the position $(1,1)$ of $A$.
In the following table in each position we write $j$ if we reach that position after having applied $S_{R,C}$ to $(1,1)$ exactly $j$ times. 		
		\begin{center}
\begin{footnotesize}
$\begin{array}{r|r|r|r|r|r|r|r|r|r|r|r|r|r|r|}&\downarrow &\downarrow &\downarrow &\downarrow &\downarrow &\downarrow &\downarrow &\downarrow &\downarrow &\downarrow &\downarrow &\uparrow &\uparrow &\uparrow\\ \hline
\rightarrow&
0 &  43 &  16 &  59  &  32 & 5 & 48 & 21 & 64 & 37 & 10 & 39 & 68 & 27\\ \hline
\rightarrow&28 & 1 & 44 & 17 & 60 & 33 & 6 & 49 & 22 & 65 & 38 & 67 & 26 & 55\\ \hline
\rightarrow&56  & 29 &  2 & 45 & 18 & 61 & 34 &  7 & 50 & 23 & 66 & 25 & 54 & 13\\ \hline
\rightarrow&14 & 57 & 30 & 3 & 46 & 19 & 62 & 35 & 8 & 51 & 24 & 53 & 12 & 41\\ \hline
\rightarrow&42 & 15 & 58 & 31 & 4 & 47 & 20 & 63 & 36 & 9 & 52 & 11 & 40 & 69\\ \hline
    \end{array}$
		\end{footnotesize}
				\end{center}
It is easy to see that we obtain $L(1,1)=(\Gamma_1, \Gamma_4, \Gamma_2, \Gamma_5, \Gamma_3)$.
  \end{ex}

\begin{thm}\label{rettangolopieno}
Let us consider a totally filled $n\times m$ array $A$. Then there exists a solution of $P(A)$ if and only if $n$ and $m$ are not both even.
\end{thm}
\proof
If $m$ and $n$ are both even, the second necessary condition of Theorem  \ref{prop:necc} is not satisfied, hence $P(A)$ has no solution.
So suppose that $m$ and $n$ are not both even.

Firstly suppose $m$ odd. Let $l$ be such that $m-2l=1$. Because of Lemma \ref{Coprimi} we have that $R:=(1,1,\dots,1)$ and $C:=(1,\dots,c_{m-l}=1,c_{m-l+1}=-1,\dots,-1)$ are a solution of $P(A)$.

Suppose now $m$ even, which implies that $n$ is odd. Let $l$ be such that $m-2l=2$.
Because of Lemma \ref{Coprimi} we have that $R:=(1,1,\dots,1)$ and $C:=(1,\dots,c_{m-l}=1,c_{m-l+1}=-1,\dots,-1)$ are a solution of $P(A)$.
\endproof

\section{Square arrays with $k$ filled diagonals}\label{SecDiag}
In this section we focus  on partially filled square arrays having the same number of filled cells in each row and column.
Let $A$ be a square array of size $n$, for $i=1,\ldots,n$ we define the \emph{$i$-th diagonal} $D_i$ to be the set of cells
$$D_i:=\{(i,1),(i+1,2),\ldots,(i-1,n)\}$$
where all the arithmetic is performed in $\mathbb{Z}_n$ using the reduced residues $\{1,2,\ldots,n\}$.
We will say that $D_i,D_{i+1},\ldots, D_{i+k-1}$ are $k$ \emph{consecutive diagonals}.

\begin{defin}
Let $k\geq1$ be an integer. We will say that
a square array $A$ of size $n\geq k$ is \emph{$k$-diagonal} if the non empty cells of $A$ are exactly those of $k$ diagonals $D_{i_1}, D_{i_2},\ldots, D_{i_k}$.
\end{defin}
Obviously, if $n=k$, $A$ is a totally filled array which has already considered in previous section.
So, let $A$ be a $k$-diagonal array of size $n>k$. In order to find solutions of $P(A)$ we assume that $R=(1,\ldots,1)$, even if this is not a necessary condition
as shown in the following example.
\begin{ex}
\label{ex5}
Let $A$ be a $3$-diagonal array of size $7$ whose filled diagonals are $D_1,D_2,D_3$.
We note that $R:=(-1, 1, 1, 1, 1, 1, 1 )$ and $C:=(-1, -1, 1, 1, 1, -1, -1 )$ are a solution of $P(A)$ even though that $R$ is not $(1,\ldots,1)$.
In the following figure we show this graphically by labeling the filled cells as done in previous examples.

\begin{center}
$A:=\begin{array}{r|r|r|r|r|r|r|r|}
&\uparrow &\uparrow     & \downarrow   & \downarrow   & \downarrow   & \uparrow   & \uparrow   \\ \hline
\leftarrow&0 &    &    &    &    &  15 &  18 \\\hline
\rightarrow &4 &  17 &       &    &   &  &  20 \\\hline
\rightarrow &16 &  11 &  3 &     &    &  &  \\\hline
\rightarrow & &  5 &  12 &  10 &       &    &   \\\hline
\rightarrow  &   &    &   6 &  13 &   2 &    &    \\\hline
\rightarrow   &  &     &   &  7 &  14 &  9 &    \\\hline
\rightarrow  & &      &   &  &  8 &  19 &  1 \\\hline
    \end{array}$
\end{center}
\end{ex}
The motivation of this choice for $R$ is explained in the following remark that we state in the case of square arrays, but it holds also
for the rectangular ones.

\begin{rem}
Let $A$ be an array of size $n$ such that the cells $(i,j),(i,j+1),(i+1,j),(i+1,j+1)$ are filled
$$\begin{array}{c|c|c|c|c|}
 &\dots & \hspace{0.28cm} j {\hspace{0.28cm}}& j+1  & \dots \\ \hline
\vdots&\ddots &  \vdots & \vdots  &\ddots  \\ \hline
i&\cdots&  \bullet & \bullet  & \cdots \\ \hline
 i+1&\cdots&\bullet   &\bullet   &\cdots  \\ \hline
\vdots&\ddots &  \vdots & \vdots  &\ddots  \\ \hline
    \end{array}$$

 Then we have that
  \begin{itemize}
    \item   $R:=(r_1,\ldots,r_{i-1},r_i=1,r_{i+1}=-1,r_{i+2},\ldots,r_n)$
  and
   $C:=(c_1,\ldots,c_{j-1},c_j=-1,c_{j+1}=1,c_{j+2},\ldots,c_n)$
  are not a solution of $P(A)$ since $S^2_{R,C}((i,j))=(i,j)$;
    \item $R:=(r_1,\ldots,r_{i-1},r_i=-1,r_{i+1}=1,r_{i+2},\ldots,r_n)$
  and
   $C:=(c_1,\ldots,c_{j-1},c_j=1,c_{j+1}=-1,c_{j+2},\ldots,c_n)$
  are not a solution of $P(A)$ since $S^2_{R,C}((i+1,j))=(i+1,j)$.
  \end{itemize}
\end{rem}

So, sign changes in both rows and columns can create very short paths on the board whenever there are a relevant number of adjacent filled positions. For this reason, from now on we look for solutions with $R=(1,\ldots,1)$. Obviously, since $A$ is a $k$-diagonal array of size $n>k$, this implies that
$C\neq(1,\ldots,1)$. We define $E=(e_1,\ldots,e_t)$, where $e_1<e_2<\ldots <e_t$, to be the list of the positions of $-1$'s in $C$.
For instance if $C=(1,-1,1,1,-1,-1)$ we have $E=(2,5,6)$.

\subsection{Cyclically $k$-diagonal square arrays}
In this subsection we investigate the case in which the $k$ filled diagonals $D_{i_1}, D_{i_2},\ldots, D_{i_k}$ of $A$ are consecutive.

\begin{defin}\label{def:kcyclic}
Let $k$ be an integer.
A square array $A$ of size $n> k$  is said to be \emph{cyclically $k$-diagonal} if
it is $k$-diagonal and the non empty diagonals are consecutive (modulo $n$).
\end{defin}

\begin{ex}\label{ex6}
The following is a cyclically $5$-diagonal array of size $9$ whose filled diagonals are $D_8,D_9, D_1,D_2,D_3$.
\begin{center}
$\begin{array}{|r|r|r|r|r|r|r|r|r|}\hline
\bullet &  \bullet &  \bullet &    &   &    &    &  \bullet &  \bullet \\ \hline
    \bullet &  \bullet &  \bullet &  \bullet &    &    &    &    &  \bullet \\ \hline
     \bullet &   \bullet &   \bullet &   \bullet &   \bullet &    &    &    &    \\ \hline
      &   \bullet &   \bullet &   \bullet &   \bullet &  \bullet &    &    &    \\ \hline
      &    &  \bullet &  \bullet &  \bullet &  \bullet &  \bullet &    &    \\ \hline
      &    &    &  \bullet &  \bullet &  \bullet &  \bullet &  \bullet &    \\ \hline
      &    &    &    &  \bullet &  \bullet &  \bullet &  \bullet &  \bullet \\ \hline
    \bullet &    &    &    &    &  \bullet &  \bullet &  \bullet &  \bullet \\ \hline
    \bullet &  \bullet &    &    &    &    &  \bullet &  \bullet &  \bullet \\ \hline
    \end{array}
    $
\end{center}
  \end{ex}

Given  a cyclically $k$-diagonal array $A$, we will say that $A$ is written in the standard form if its nonempty diagonals are $D_1,\dots,D_k$.
Note that the array of Example \ref{ex5} is in standard form, while the array in Example \ref{ex6} is not.
In the following, by Remark \ref{11=ij}, we can assume without loss of generality that a cyclically $k$-diagonal array $A$ is written in the standard form.

We start, determining a necessary condition for the existence of a solution.
\begin{prop}\label{nkodd}
Let $A$ be a cyclically $k$-diagonal
array of size $n> k$.
If $P(A)$ admits a solution then $n$ and $k$ are both odd and $k\neq1$.
\end{prop}
\proof
Since $|F(A)|=kn$, from Corollary \ref{cor:necc} we have $n$ and $k$ odd.
Also, note that if $k=1$, then every cell of the filled diagonal of $A$ is a closed subarray. Therefore the thesis follows from Theorem \ref{prop:necc}.
\endproof
Hence, in the following let $k\geq3$ be an odd integer.

\begin{lem}\label{Percorso1}
Let $k\geq3$ be an odd integer and let $A$ be a cyclically $k$-diagonal array of size $n> k$.
Then the vectors $R:=(1,\dots,1)$ and $C\in \{-1,1\}^n$, whose $-1$ are in $E$, are a solution of $P(A)$ if and only if:
\begin{itemize}
\item[1)] denoted by $d=\gcd(n,k-1)$, the list $E$ covers all the congruence classes modulo $d$;
\item[2)] the list $L(1,1)$ covers all the positions of $\{(e,e)\ |\ e\in E \}$.
\end{itemize}
\end{lem}
\proof
If $R$ and $C$ are a solution, the list $L(1,1)$ covers all the positions of $F(A)$ and hence all the positions of $\{(e,e)\ |\ e\in E \}\subseteq D_1$.
Let us suppose, by absurd, that $R$ and $C$ are a solution and that there exists a congruence class $[f]$ modulo $d$ that is not covered by $E$.
We note that, for $(i,i)\in D_1$ and $i-(k-1)\not\in E$, we have $S_{R,C}((i,i))=(i-(k-1),i-(k-1))\in D_1$ where $i-(k-1)\equiv i \pmod{d}$.
Therefore the list $L(f,f)$ is contained in $D_1$, which is absurd.

Conversely, let us suppose, by absurd, that
the hypotheses 1) and 2) of the statement are satisfied, but
 $L(1,1)$
 does not cover $F(A)$. Let $(f,g)\in F(A)\setminus L(1,1)$ and consider the list  $L(f,g)$.
Obviously $L(1,1)\cap L(f,g)=\emptyset$,
hence,  by hypothesis 2)
\begin{equation}\label{contraLemma2}
\{(e,e)| e\in E\}\cap L(f,g)=\emptyset.
\tag{$\ast$}
  \end{equation}
We note that, given a position $(i,j)$, $S_{R,C}((i,j))\in D_{1}$ only in the following two cases:
\begin{itemize}
\item[A)] $(i,j)\in D_1$ and $i\not\in \{e+(k-1)\ |\ e\in E\}$;
\item[B)] $(i,j)\in D_{3}$ and $S_{R,C}((i,j))=(e,e)$ for some $e\in E$.
\end{itemize}

Since from A) for $(i,i)\in D_1\setminus \{(e+(k-1),e+(k-1))\ |\ e\in E\}$ we have $S_{R,C}((i,i))=(i-(k-1),i-(k-1))$ and, by hypothesis 1),
the list $E$ covers all the congruence classes modulo $d$,
it follows that $L(f,g)\not \subseteq D_1$.
Suppose now that there exists $(i',j')\in L(f,g)\setminus D_1$ such that $S_{R,C}((i',j'))\in D_1$.
From $B)$ it follows that  $S_{R,C}((i',j'))\in \{(e,e)\ |\ e\in E\}$, but obviously $S_{R,C}((i',j'))\in L(f,g)$, which is in contradiction
with \eqref{contraLemma2}.
Hence $L(f,g)\cap D_1=\emptyset$.

Now we note that, given $(i,j)\in D_h$ with $h\not= 1$, we have $S_{R,C}((i,j))=(i+1,j+1)\in D_h$ if $j+1\not\in E$ and
$S_{R,C}((i,j))\in D_{h-2}$ if $j+1\in E$ where the subscripts of the diagonals are considered modulo $k$. Since $k$ is odd, this means that the list $L(f,g)$ will reach every non empty diagonal and hence also the diagonal $D_1$,
but this is absurd because $L(f,g)\cap D_1=\emptyset$.
\endproof
As a consequence we have the following result.

\begin{prop}\label{corCMPP}
Let $k\geq3$ be an odd integer and let
$A$ be a cyclically $k$-diagonal
array of size $n> k$.
If $\gcd(n,k-1)=1$, then the vectors $R:=(1,\dots,1)$ and $C:=(-1,1\dots,1)$ are a solution of $P(A)$.
\end{prop}
We point out that the result of Proposition \ref{corCMPP} was previously obtained in \cite{CMPPGlobally} (see Proposition 3.4).

Given a cyclically $k$-diagonal array of size $n> k$ and vectors $R=(1,\dots,1)$ and $C\in \{-1,1\}^n$, whose $-1$ are in positions $E=(e_1,\dots,e_t)$
where $e_1<e_2<\dots<e_t$, we would like to study some properties of the move function $S_{R,C}$.
In particular, considered an element $(e,e)\in D_1$ with $e\in E$, there exists a minimum $m\geq 1$ such that $S_{R,C}^m((e,e))=(e',e')$ for some $e'\in E$.
We define the permutation $\omega_C$ on $E$ as $\omega_C(e)=e'$.
Similarly, given $e\in E$, there exists a minimum $m\geq 1$ such that $e-m(k-1)\equiv e''\pmod{n}$ for some $e''\in E$.
We define the permutation $\omega_{1,C}$ on $E$ as $\omega_{1,C}(e)=e''$.
Finally we define the permutation $\omega_{2,C}$ on $E=(e_1,\ldots,e_t)$ as $\omega_{2,C}(e_i)=e_{i+(k-1)}$ where the indices are considered modulo $t$.

\begin{lem}\label{Percorso2}
Let $k\geq3$ be an odd integer and let $A$ be a cyclically $k$-diagonal array of size $n> k$.
Then the vectors $R:=(1,\dots,1)$ and $C\in \{ -1,1\}^n$, whose $-1$ are in $E$, are a solution of $P(A)$ if and only if:
\begin{itemize}
\item[1)] denoted by $d=\gcd(n,k-1)$, the list $E$ covers all the congruence classes modulo $d$;
\item[2)] the permutation $\omega_{2,C}\circ \omega_{1,C}$ on $E$ is a cycle of length $t=|E|$.
\end{itemize}
\end{lem}
\proof
We note that $\omega_C=\omega_{2,C}\circ \omega_{1,C}$ and that $\omega_C$ is a cycle of length $t$ if and only if the list
$L(1,1)$ covers all the positions of $\{(e,e)\ |\ e\in E\}$. Then the claim follows from Lemma \ref{Percorso1}.
\endproof

\begin{prop}\label{extension}
Let $k\geq3$ be an odd integer and let
$A$ be a cyclically $k$-diagonal
array of size $n> k$. Let us assume that the vectors $R:=(1,\dots,1)$ and $C$ are a solution of $P(A)$.
Then there exists a solution of $P(A')$ for any cyclically $k$-diagonal array $A'$ of size $n'=n+\lambda(k-1)$ for any integer $\lambda\geq 0$.
\end{prop}
\proof
Let us consider the vector $R':=(1,\dots,1)$ and the vector $C'\in \{-1,1\}^{n'}$ that has the $-1$ in the same positions, denoted by $E=(e_1,\dots,e_t)$, of $C$.
Then the permutations $\omega_{1,C}$ and $\omega_{1,C'}$ on $E$ are identical and the same holds for the permutations $\omega_{2,C}$ and $\omega_{2,C'}$.
Thus $\omega_C=\omega_{C'}$. Since $\gcd(n,k-1)=\gcd(n+\lambda(k-1),k-1)$, the claim follows from Lemma \ref{Percorso2}.
\endproof
As a consequence of Proposition \ref{extension} we obtain the following result:

\begin{thm}\label{smallcases}
Let $3\leq k< 200$ be an integer and
let $A$ be a cyclically $k$-diagonal
array of size $n> k$. Then there exists a solution of $P(A)$ if and only if $n$ and $k$ are both odd.
\end{thm}
\proof
The necessary condition $n$ and $k$ odd follows from Proposition \ref{nkodd}.
For any odd integer $k$, with $3\leq k< 200$, we checked, using a computer, that the odd values $n\in [k+2,2k-1]$ have a solution with vector $R=(1,\dots,1)$.
Hence the thesis follows from Proposition \ref{extension}.
\endproof

At a first reading it seems that, in Lemma \ref{Percorso2} (and in the analogous Lemma \ref{Percorso1}), we should also require that $n$ is not an even number,
according to Proposition  \ref{nkodd}. In the next result we show that this condition is hidden in the hypothesis of those lemmas.
\begin{prop}
Let $k\geq3$ be an odd integer and let $A$ be a cyclically $k$-diagonal array of size $n>k$.
Let us consider the vectors $R:=(1,\ldots,1)$ and $C\in \{-1,1\}^n$, whose $-1$ are in $E$.
Let us suppose that, denoted by $d=\gcd(n,k-1)$, the list $E$ covers all the congruence classes modulo $d$.
Then the permutation $\omega_{2,C}\circ \omega_{1,C}$ on $E$ can be a cycle of length $t=|E|$ only if $n$ is odd.
\end{prop}
\proof
Suppose, by absurd, that $n$ is even (which implies $d$ even) and that $\omega_C=\omega_{2,C}\circ \omega_{1,C}$ is a cycle of length $t$.
Denoted by $\tau$ the permutation of the list $E$ such that $\tau(e_i)=e_{i+1}$, we have $\omega_{2,C}=\tau^{k-1}$. Therefore $\omega_{2,C}$ is always an even permutation.
Let us partition $E$ in the congruence classes modulo $d$ by setting $E^h:=\{e_i\in E| e_i\equiv h \pmod{d} \}.$
Since $d=\gcd(n,k-1)$, the permutation $\omega_{1,C}$ is the product of cyclic permutations on the sets $E^h$, each of which has length $|E^h|$.
Since $d$ is even and no $E^h$ is empty, $\omega_{1,C}$ has parity $\sum_{h=0}^{d-1} (|E^h|-1)\equiv (\sum_{h=0}^{d-1}|E^h|)-d\equiv t \pmod{2}$, which is also the parity of $\omega_{2,C}\circ \omega_{1,C}$.
By hypothesis we have assumed that $\omega_C$ is a cycle of length $t$,
then it has parity $t-1$, but this is absurd because $\omega_C=\omega_{2,C}\circ \omega_{1,C}$ should have parity $t$.
\endproof
\begin{rem}
  Since the hypotheses of Lemma $\ref{Percorso1}$ imply those of Lemma $\ref{Percorso2}$,
  the condition $n$ odd is also hidden in Lemma $\ref{Percorso1}$.

\end{rem}
Now we introduce a more general application of Lemma \ref{Percorso2} that helps us in finding a solution of the \probname\ whenever $n$ is sufficiently large.
\begin{thm}
Let $k\geq3$ be an odd integer and let $A$ be a cyclically $k$-diagonal
array of odd size $n\geq (k-2)(k-1)$. Then there exists a solution of $P(A)$.
\end{thm}
\proof
Set again $R:=(1,\dots,1)$ and let $d=\gcd(n,k-1)$, because of Proposition \ref{corCMPP} we can assume that $d\geq 3$.
Now we want to define a vector $C$ that satisfies the hypothesis of Lemma \ref{Percorso2}.
  Since $n\geq (k-2)(k-1)\geq (k-d+1)(k-1)$, there exist $k$ columns $e_1,\dots,e_k$ such that:
\begin{itemize}
\item $e_i= i $ for any $i\in[1,d-1]$;
\item $e_{i}=d+(k-1)(i-d)$ for any $i \in [d, k]$.
\end{itemize}
Consider now the vector $C\in \{-1,1\}^n$, whose $-1$ are in positions $(e_1,\dots,e_k)$. Clearly the set $\{e_1,\dots,e_k\}$ covers all the congruence classes modulo $d$.
Since the $e_i$'s in the congruence class of $d$ are at distance $k-1$, $\omega_{1,C}$ is the cyclic permutation $(e_{d},e_{k},e_{k-1},\dots, e_{d+1})$.
In this context $t=|E|=k$ and $k-1\equiv -1\pmod{t}$, hence we have that $\omega_{2,C}$ is the cyclic permutation $(e_{1},e_{k},e_{k-1},\dots, e_{2})$.
Since $k-d+1$ is odd, one can easily check that $\omega_{2,C}\circ \omega_{1,C}$ is given by the cyclic permutation:
$$(e_{1},e_{k},e_{k-2},\dots,e_{d+2}, e_d,e_{k-1},e_{k-3},\dots, e_{d+1},e_{d-1},\dots, e_{2}).$$
Therefore, because of Lemma \ref{Percorso2}, the vectors $R$ and $C$ are a solution for the \probname\ for any cyclically $k$-diagonal array of size $n\geq (k-2)(k-1)$.
\endproof

\subsection{$k$-diagonal square arrays}
In this subsection we consider the case in which the filled diagonals of $A$ are not necessarily consecutive.
In order to study the solutions of $P(A)$ it is important to know the number of empty diagonals between two filled diagonals,
 we thus introduce the notion of \emph{empty strip}.

\begin{ex}\label{strisce}
Here we have a $5$-diagonal array of size $11$, the non empty diagonals are $D_1,D_4,\\D_6,D_7$ and $D_{11}$.
$$\begin{array}{|r|r|r|r|r|r|r|r|r|r|r|r|}\hline
     \bullet &   \bullet &   & &  & \bullet   &  \bullet &   & \bullet &    &   \\\hline
   &   \bullet & \bullet &  & &  & \bullet   &  \bullet &   & \bullet &    \\\hline
	 &  &   \bullet & \bullet &  & &  & \bullet   &  \bullet &   & \bullet   \\\hline
	\bullet & &  &   \bullet & \bullet &  & &  & \bullet   &  \bullet &      \\\hline
	&	\bullet & &  &   \bullet & \bullet &  & &  & \bullet   &  \bullet      \\\hline
	\bullet &	&	\bullet & &  &   \bullet & \bullet &  & &  & \bullet       \\\hline
	\bullet &	\bullet &	&	\bullet & &  &   \bullet & \bullet &  & &       \\\hline
		& \bullet &	\bullet &	&	\bullet & &  &   \bullet & \bullet &  &     \\\hline
		&	& \bullet &	\bullet &	&	\bullet & &  &   \bullet & \bullet &      \\\hline
		&	&	& \bullet &	\bullet &	&	\bullet & &  &   \bullet & \bullet    \\\hline
		\bullet &	&	&	& \bullet &	\bullet &	&	\bullet & &  &   \bullet   \\\hline
			\end{array}$$
\end{ex}
Let $A$ be a $k$-diagonal array of size $n> k$. A set $S=\{D_{r+1},D_{r+2},\ldots, D_{r+t}\}$ is said to be
an \emph{empty strip of width $t$} if  $D_{r+1},D_{r+2},\ldots,D_{r+t}$ are empty diagonals, while $D_r$ and $D_{r+t+1}$ are filled diagonals.
The array of Example \ref{strisce} has three empty strips $S_1=\{D_2,D_3\}$, $S_2=\{D_5\}$ and $S_3=\{D_8,D_9,D_{10}\}$.
\begin{defin}
Let $A$ be a $k$-diagonal array of size $n > k$.
We will say that $A$ is a \emph{$k$-diagonal array with width $s$}
if all the empty strips of $A$ have width $s$.
\end{defin}
\begin{ex}\label{strisceuguali}
The following is a $5$-diagonal array of size $11$ with width $3$.
$$\begin{array}{|r|r|r|r|r|r|r|r|r|r|r|r|}\hline
     \bullet &    &   & & \bullet & \bullet   &  &   &  &  \bullet  & \bullet  \\\hline
    \bullet & \bullet &    &   & & \bullet & \bullet   &  &   &  &  \bullet   \\\hline
	\bullet &	 \bullet & \bullet &    &   & & \bullet & \bullet   &  &   &   \\\hline
&	\bullet &	 \bullet & \bullet &    &   & & \bullet & \bullet   &  &    \\\hline
& &	\bullet &	 \bullet & \bullet &    &   & & \bullet & \bullet   &   \\\hline
& & &	\bullet &	 \bullet & \bullet &    &   & & \bullet & \bullet  \\\hline
\bullet & & & &	\bullet &	 \bullet & \bullet &    &   & & \bullet  \\\hline
\bullet &\bullet & & & &	\bullet &	 \bullet & \bullet &    &   &   \\\hline
& \bullet &\bullet & & & &	\bullet &	 \bullet & \bullet &    &  \\\hline
& & \bullet &\bullet & & & &	\bullet &	 \bullet & \bullet &  \\\hline
& & &  \bullet &\bullet & & & &	\bullet &	 \bullet & \bullet \\\hline
			\end{array}$$
\end{ex}
We will say that a $k$-diagonal array $A$ is written in standard form if $D_1$
is a filled diagonal and $D_n$ is an empty diagonal.
For instance the array of Example \ref{strisce} is not in a standard form, while that of Example
\ref{strisceuguali} is written in a standard form.
In the following we will suppose
that $A$ is written in a standard form, since this is not restrictive in order to study the solution of $P(A)$.

Now we present a result which generalizes Lemma \ref{Percorso1}.
\begin{lem}\label{Percorso1bis}
Let $k\geq3$ be an odd integer and let $A$ be a $k$-diagonal array of size $n> k$ and width $s$.
Then the vectors $R:=(1,\dots,1)$ and $C\in \{-1,1\}^n$, whose $-1$ are in $E$, are a solution of $P(A)$ if and only if:
\begin{itemize}
\item[1)] denoted by $d=\gcd(n,s+1)$, the list $E$ covers all the congruence classes modulo $d$;
\item[2)] the list $L(1,1)$ covers all the positions of $\{(e,e)| e\in E\}$.
\end{itemize}
\end{lem}
\proof
Let $D_{i_1}=D_1, D_{i_2},\ldots,D_{i_k}$ with $i_j<i_{j+1}$, for any $j=1,\ldots,k-1$, be the $k$ non empty diagonals of $A$.

If $R$ and $C$ are a solution, the list $L(1,1)$ covers all the filled positions of $A$ and hence all the positions of $\{(e,e)| e \in E\}\subseteq D_1$.
Let us suppose, by absurd, that $R$ and $C$ are a solution and that there exists a congruence class $[f]$ modulo $d$ that is not covered by the list $E$.
We note that, for $(i,i)\in D_1$ and $i+s+1 \not\in E$, we have $S_{R,C}((i,i))=(i+s+1,i+s+1)\in D_1$ where $i+s+1\equiv i \pmod{d}$.
Therefore the list $L(f,f)$ is contained in $D_1$, which is absurd.

Conversely, let us suppose, by absurd,  that
the hypotheses 1) and 2) of the statement are satisfied, but
 $L(1,1)$ does not cover $F(A)$.
Let $(f,g)\in F(A)\setminus L(1,1)$ and consider the list  $L(f,g)$; obviously $L(1,1)\cap L(f,g)=\emptyset$, hence by 2) we have that
 \begin{equation}\label{eqn}
\{(e,e)| e\in E\}\cap L(f,g)=\emptyset.
\tag{$\ast$}
  \end{equation}
We note that, given a position $(i,j)$, $S_{R,C}((i,j))\in D_{1}$ only in the following two cases:
\begin{itemize}
\item[A)] $(i,j)\in D_1$ and $i+s+1 \not\in E$;
\item[B)] $(i,j)\in D_{i_{3}}$ and $S_{R,C}((i,j))=(e,e)$ for some $e\in E$.
\end{itemize}
Since for $i+s+1\notin E$ we have $S_{R,C}((i,i))=(i+s+1,i+s+1)$ and, by hypothesis 1), the list $E$ covers all the congruence classes modulo $d$,
it follows that $L(f,g)\not \subseteq D_1$.
Suppose now that there exists $(i',j')\in L(f,g)\setminus D_1$ such that $S_{R,C}((i',j'))\in D_1$.
From $B)$ it follows that  $S_{R,C}((i',j'))\in \{(e,e)| e\in E\}$, but obviously $S_{R,C}((i',j'))\in L(f,g)$, which is in contradiction with \eqref{eqn}.
Hence $L(f,g)\cap D_1=\emptyset$.
Given $(i,j)\in D_{i_h}$, set $S_{R,C}((i,j))=(i'',j'')$. We note that $S_{R,C}((i,j))\in D_{i_h}$ if $j''\not\in E$, otherwise
$S_{R,C}((i,j))\in D_{i_{h-2}}$ where the subscripts of the diagonals are considered modulo $k$. Since $k$ is odd it means that the list $L(f,g)$
will reach every non empty diagonal and hence also the diagonal $D_1$, but this is absurd because  $L(f,g)\cap D_1=\emptyset$.
\endproof

The following result is an easy consequence of the previous lemma.
\begin{prop}\label{s+1,n}
Let $k\geq3$ be an odd integer and let $A$ be a $k$-diagonal array of size $n > k$ and width $s$.
If $\gcd(n,s+1)=1$, then the vectors $R:=(1,\dots,1)$ and $C:=(-1,1,\dots,1)$ are a solution of $P(A)$.
\end{prop}
\proof
It is immediate to see that conditions 1) and 2) of Lemma \ref{Percorso1bis} are satisfied.
\endproof
If the array $A$ has exactly two strips of the same width the following holds.
\begin{cor}\label{2strisce}
Let $k\geq3$ be an odd integer and let $A$ be a $k$-diagonal array of size $n > k$ and width $s=\frac{n-k}{2}$.
If $\gcd(n,k-2)=1$, then the vectors $R:=(1,\dots,1)$ and $C:=(-1,1,\dots,1)$ are a solution of $P(A)$.
\end{cor}
\proof
Since $n$ and $k$ are odd, $\gcd(n,k-2)=1$ if and only if $\gcd(n,\frac{n-k+2}{2})=1$.
Note that  $\frac{n-k+2}{2}=s+1$, hence the thesis follows from Proposition \ref{s+1,n}.
\endproof

\begin{ex}
The array of Example $\ref{strisceuguali}$ satisfies the hypotheses of Corollary \ref{2strisce} with $n=11$, $k=5$ and $s=3$.
We show the tour obtained with the solution presented in the corollary writing in a position $j$ if we reach that position after having applied $S_{R,C}$ to $(1,1)$ exactly $j$ times. 		

$$\begin{array}{r|r|r|r|r|r|r|r|r|r|r|r|r|}&  \uparrow & \downarrow   &\downarrow   & \downarrow& \downarrow & \downarrow & \downarrow &\downarrow   &\downarrow  &\downarrow   &\downarrow \\\hline
\rightarrow   &  0 &    &   & & 37 & 15  &  &   &  &  53 & 32  \\\hline
\rightarrow   & 22 & 3 &    &   & & 38 & 18  &  &   &  &  54   \\\hline
\rightarrow&	44 &	23 & 6 &    &   & & 39 & 21   &  &   &   \\\hline
\rightarrow&&	45 &	24 & 9 &    &   & & 40 & 13   &  &    \\\hline
\rightarrow&& &	46 &	 25 & 1 &    &   & & 41 & 16 &   \\\hline
\rightarrow&& & &	47 &	26 & 4 &    &   & & 42 & 19  \\\hline
\rightarrow&11 & & & &	48 & 27 & 7 &    &   & & 43  \\\hline
\rightarrow&33 & 14 & & & &	49 &	28 & 10 &    &   &   \\\hline
\rightarrow&& 34 & 17 & & & &	50 &	29 & 2 &    &  \\\hline
\rightarrow&& & 35 & 20 & & & &	51 &	30 & 5 &  \\\hline
\rightarrow&& & & 36 & 12 & & & &	52 &	31 & 8 \\\hline
			\end{array}$$
\end{ex}

\section{Cyclically almost $k$-diagonal square arrays}\label{SecAlmost}
We have seen that the necessary conditions for the existence of a solution of the problem when we have a $k$-diagonal
array of size $n$ are $n$ and $k$ odd.
Hence if we want to obtain some results about arrays of even size and/or with an even number of filled diagonals, we have to change the structure of filled cells.
In order to break the symmetry as little as possible we consider an array with $k$ filled diagonals and another ``extra'' filled position.

\begin{defin}
Let $k\geq1$ be an integer. We will say that a square array $A$ of size $n> k$ is \emph{almost $k$-diagonal}
 if the non empty cells of $A$ are exactly those of $k$ diagonals $D_{i_1}, D_{i_2},\ldots, D_{i_k}$ together with
another extra filled position belonging to a diagonal $D_j$ with $j\neq i_1,\ldots,i_k$.
\end{defin}
In other words, an almost $k$-diagonal array can be obtained by adding an extra filled cell to a $k$-diagonal array, hence it has exactly $kn+1$ filled positions.

\begin{ex}\label{almostnonstandard}
The following is an almost $2$-diagonal array of size $5$, whose totally filled diagonals are $D_3$ and $D_4$, and the extra filled
position is $(4,4)$.
$$\begin{array}{|r|r|r|r|r|}\hline
   &    &  \bullet & \bullet   &   \\\hline
 &    &  & \bullet & \bullet   \\\hline
 \bullet &    &  & & \bullet     \\\hline
   \bullet  &  \bullet & & \bullet   &   \\\hline
       &  \bullet & \bullet   &   & \\\hline	
		\end{array}$$
\end{ex}

First of all, we determine a necessary condition for the existence of a solution of $P(A)$ when $A$ is an almost $k$-diagonal array.
\begin{prop}\label{nkeven}
Let $A$ be an almost $k$-diagonal
array of size $n> k$.
If $P(A)$ admits a solution then $kn$ is even with $k\neq1$ or $(n,k)=(2,1)$.
\end{prop}
\proof
Since $|F(A)|=kn+1$, from Corollary \ref{cor:necc} we have that $kn$ must be even.
Then note that for $k=1$, $A$ admits nontrivial closed subarrays except when $n=2$.
Hence the thesis follows from Theorem \ref{prop:necc}.
\endproof

\begin{defin}
  Let $A$ be an almost $k$-diagonal array of size $n>k$. We will say that $A$ has \emph{width} $s$ if $A$
can be obtained adding an extra filled cell to a $k$-diagonal array with width $s$.
\end{defin}

An almost $k$-diagonal array $A$ of width $s$ is said to be in standard form if the extra filled cell is $(1,\ell)$ with $2\leq \ell\leq s+1$,
$D_1$ is a totally filled diagonal,
$D_n$ has at most $(1,2)$ as filled position (note that this happens only when $\ell=2$).
As before, in order to study the solution of the problem it is not restrictive to suppose that $A$ is written in standard form.

Also in this section we will look for solutions with $R=(1,\ldots,1)$. As done in Section \ref{SecDiag}, by $E$ we will denote the
list of positions of $-1$'s in $C$. Moreover by $\mathcal{C}_h$ we will mean the set of the cells $(j,j)$ of $D_1$ such that $j\equiv h \pmod d$,
where  $d=\gcd(n,s+1)$.

\begin{lem}\label{nevenkodd}
Let $k\geq3$ be an odd integer and let $A$ be an almost $k$-diagonal array of size $n> k$ and with width $s$, whose extra filled position is $(1,\ell)$.
Then the vectors $R:=(1,\dots,1)$ and $C\in \{-1,1\}^n$, whose $-1$ are in $E$ and with $c_\ell=1$,
are a solution of $P(A)$ if and only if:
\begin{itemize}
\item[1)] denoted by $d=\gcd(n,s+1)$, the list $E$ covers all the congruence classes modulo $d$ if $\ell\equiv1\pmod d$,
otherwise the list $E$ covers all the congruence classes modulo $d$ except at most one among the classes of $1$ and  of $\ell$;
\item[2)] the list $L(1,1)$ covers all the positions of $\{(e,e)| e\in E\}$.
\end{itemize}
\end{lem}
\proof
Let $D_{i_1}=D_1, D_{i_2},\ldots,D_{i_k}$ with $i_j<i_{j+1}$, for any $j=1,\ldots,k-1$, be the $k$ totally filled diagonals of $A$.

If $R$ and $C$ are a solution, the list $L(1,1)$ covers all the filled positions of $A$ and hence all the positions of $\{(e,e)| e\in E\}\subseteq D_1$.
Let us suppose, by absurd, that $R$ and $C$ are a solution and that one of the following two conditions holds:
\begin{itemize}
\item[a)] there exists a congruence class $[f]$ modulo $d$, with $f\not\equiv 1,\ell \pmod d$, that is not covered by the list $E$;
\item[b)] the congruence classes of $1$ and of  $\ell$ modulo $d$ are not covered by the list $E$.
\end{itemize}
Firstly suppose we are in case a).
We note that, for $(i,i)\in D_1$ with $i+(s+1)\not\in E$ and $i\not\equiv 1\pmod d$, we have $S_{R,C}((i,i))=(i+(s+1),i+(s+1))\in D_1$ where $i+(s+1)\equiv i \pmod{d}$.
Therefore the list $L(f,f)$ is contained in $D_1$, which is absurd.

Now suppose we are in case b).
Note that $S_{R,C}((1,1))=(\ell,\ell)\in \mathcal{C}_\ell$ and that
 if $(j,j)\in \mathcal{C}_\ell$ then $S_{R,C}((j,j))\in \mathcal{C}_\ell$ except when $j=n+\ell-(s+1)$, in fact $S_{R,C}((n+\ell-(s+1),n+\ell-(s+1))=(1,\ell)$.
Also, since $A$ is in standard form, $S_{R,C}((1,\ell))=(s+2,s+2)\in \mathcal{C}_1$ and if $(j,j)\in \mathcal{C}_1$ then $S_{R,C}((j,j))\in \mathcal{C}_1$ except when $j=1$.
Hence if $\ell\not\equiv 1\pmod d$, then $L(1,1)\subseteq \mathcal{C}_1\cup \mathcal{C}_\ell \cup (1,\ell)$, otherwise $L(1,1)\subseteq \mathcal{C}_1\cup (1,\ell)$. In both cases we obtain a contradiction.

Conversely, let us suppose, by absurd, that conditions 1) and 2) of the statement are satisfied, but
$L(1,1)$ does not cover $F(A)$, and let $(f,g)\in F(A)\setminus L(1,1)$.
Let us consider the list  $L(f,g)$, obviously
 $L(1,1)\cap L(f,g)=\emptyset$, hence by 2) we have that
 \begin{equation}\label{e}
\{(e,e)| e\in E\}\cap L(f,g)=\emptyset.
\tag{$\ast$}
  \end{equation}
We note that, given a position $(i,j)$, $S_{R,C}((i,j))\in D_{1}$ only in the following four cases:
\begin{itemize}
\item[A)] if $(i,j)=(i, i)\in D_1$ and $i\not\in \{e-(s+1)| e\in E\}\cup\{1,\ell-(s+1)\}$, then $S_{R,C}((i,i))=(i+(s+1),i+(s+1))$;
\item[B)] if $(i,j)=(1,1)$, then $S_{R,C}((1,1))=(\ell,\ell)$;
\item[C)] if $(i,j)=(1,\ell)$ and $s+2\not\in E$, then $S_{R,C}((1,\ell))=(s+2,s+2)$;
\item[D)] if $(i,j)\in D_{i_3}$ and $S_{R,C}((i,j))=(e,e)$ for some $e\in E$.
\end{itemize}
If $f\not\equiv 1,\ell\pmod d$, from A) and hypothesis 1),
it follows that $L(f,g)\not \subseteq D_1\cup\{(1,\ell)\}$.
Now let $f\equiv 1,\ell \pmod d$ and assume $L(f,g)\subseteq D_1\cup\{(1,\ell)\}$.
From A), B) and C) this implies that $L(f,g)=\mathcal{C}_1 \cup \mathcal{C}_\ell \cup \{(1,\ell)\}$.
But this happens if and only if the list $E$ does not cover  the congruence classes of $1$ and of  $\ell$ modulo $d$, which is in contradiction with hypothesis 1).
Hence, in any case, $L(f,g)\not \subseteq D_1\cup\{(1,\ell)\}$.

Suppose now that there exists $(i',j')\in L(f,g)\setminus \{D_1\cup \{(1,\ell)\} \}$ such that $S_{R,C}((i',j'))\in D_1 \cup \{(1,\ell)\}$.
We note that $S_{R,C}((i',j'))=(1,\ell)$ if and only if $(i',j')=(n+\ell-(s+1),n+\ell-(s+1))\in D_1$, but $(i',j')\not\in D_1$.
From D) it follows that  $S_{R,C}((i',j'))\in \{(e,e)| e\in E\}$, but obviously $S_{R,C}((i',j'))\in L(f,g)$, which is in  contradiction with \eqref{e}.
Hence $L(f,g)\cap (D_1\cup\{(1,\ell)\})=\emptyset$.

Set $\delta_h=s+1$ if $D_{i_h-1}$ is empty and $\delta_h=1$ otherwise. We note that, given $(i,j)\in D_{i_h}$ with $h=2,\ldots,k$, we have $S_{R,C}((i,j))=(i+\delta_h,j+\delta_h)\in D_{i_h}$ if $j+\delta_h\not\in E$ and
$S_{R,C}((i,j))\in D_{i_{h-2}}$ if $j+\delta_h\in E$ where the subscripts of the diagonals are considered modulo $k$. Since $k$ is odd it means that the list $L(f,g)$
will reach every totally filled diagonal and hence also the diagonal $D_1$, but this is absurd because $L(f,g)\cap \{D_1\cup \{(1,\ell)\}\}=\emptyset$.
\endproof
As an immediate consequence we have the following result.

\begin{prop}
  Let $k\geq3$ be an odd integer and let
$A$ be an almost $k$-diagonal
array of size $n> k$ and with width $s$, whose extra filled position is $(1,\ell)$ with $\ell$ even.
If $\gcd(n,s+1)=2$, then the vectors $R:=(1,\dots,1)$ and $C:=(-1,1,\dots,1)$ are a solution of $P(A)$.
\end{prop}

We have an analogous partial result in the case $\ell$ odd only for almost $k$-diagonals arrays in which the $k$ totally filled
diagonals are consecutive. Hence we introduce the following concept.
\begin{defin}
  An almost $k$-diagonal array $A$ whose $k$ totally filled diagonals are consecutive will be said a \emph{cyclically almost $k$-diagonal array}.
\end{defin}
Note that in this case if we write the array in standard form the extra filled position is $(1,\ell)$ with $2\leq \ell \leq n-k+1$ and that the width is $s=n-k$.

\begin{ex}\label{almoststandard}
This is a cyclically almost $2$-diagonal array of size $5$ with
 extra filled position $(1,4)$ written in standard form.
$$\begin{array}{|r|r|r|r|r|}\hline
       \bullet & & &\bullet  & \bullet \\\hline	
    \bullet & \bullet  & & &   \\\hline
   & \bullet & \bullet & &   \\\hline
     &  & \bullet & \bullet &     \\\hline
    & & &  \bullet &  \bullet    \\\hline
		\end{array}$$
\end{ex}

From now on we will consider a \underline{cyclically} almost $k$-diagonal array written in standard form.
\begin{prop}\label{3mod4}
  Let $k\equiv 3\pmod{4}$ and let
$A$ be a cyclically almost $k$-diagonal
array of size $n> k$ with the extra filled position $(1,\ell)$ with $\ell$ odd.
If $\gcd(n,k-1)=2$, then the vectors $R:=(1,\dots,1)$ and $C:=(1,-1,1,\dots,1,c_{\ell}=-1,1,\ldots,1)$ are a solution of $P(A)$.
\end{prop}
\proof
Let $R$ and $C$ be as in the statement.
For any $i=2,\ldots,k$ we set
$D_i^1=\{(x,y)\in D_i\ |\ 2\leq y\leq \ell-1\}$ and $D_i^2= D_i\setminus D_i^1$.
It is easy to see that
$$L(1,1)=((1,1),D_k^2,D_{k-2}^1,D_{k-4}^2,D_{k-6}^1,\ldots,D_3^2,(2,2),\ldots,(k,k)).$$
Set $B=((2,2),\ldots,(k,k))\subseteq L(1,1)$.
Note that $(2,2)\in D_1$ hence either one between $(k+1,k+1)$ and $(\ell+k-1,\ell+k-1)$ is in $B$, or $B\subseteq D_1$.
It is not hard to see that $B\subseteq D_1$ implies that
$B=(b_1,b_2,\ldots,b_{|B|})$ where
$b_\alpha=(2-(\alpha-1)(k-1),2-(\alpha-1)(k-1))$ for $\alpha=1,\ldots,|B|$.
But $ 2-(\alpha-1)(k-1)\neq k$ since $k$ is odd.
So  $B\not\subseteq D_1$, hence one between $(k+1,k+1)$ and $(\ell+k-1,\ell+k-1)$ is in $B$. Since $\ell+k-1$ is odd
 while $k+1$ is even, we have that $B=((2,2),\ldots,(k+1,k+1),\ldots,(k,k))$ with $((2,2),\ldots,(k+1,k+1))\subseteq D_1$ which implies:
$$L(1,1)=((1,1),D_k^2,D_{k-2}^1,D_{k-4}^2,D_{k-6}^1,\ldots,D_3^2,(2,2),\ldots,(k+1,k+1),$$
$$D_{k-1}^1,D_{k-3}^2,\ldots, D_4^2,D_2^1,(1,\ell),(n-k+2,n-k+2),\dots,(k,k)).$$
Set $G=((n-k+2,n-k+2),\ldots,(k,k))\subseteq L(1,1)$.
Note that $(n-k+2,n-k+2)\in D_1$ hence either $(\ell+k-1,\ell+k-1)\in G$, or $G\subseteq D_1$.
It is not hard to see that $G\subseteq D_1$ implies that
$G=(g_1,g_2,\ldots,g_{|G|})$ where
$g_\alpha=(n-k+2-(\alpha-1)(k-1),n-k+2-(\alpha-1)(k-1))$ for $\alpha=1,\ldots,|G|$.
But $ n-k+2-(\alpha-1)(k-1)=n+1-\alpha(k-1)\equiv k\pmod{n}$ if and only if $\alpha\equiv n/2-1\pmod{n/2}$ that is $G$ coincides with the odd elements of $D_1$ except $(1,1)$. So we can suppose that $(\ell+k-1,\ell+k-1)\in G$ which implies:
$$L(1,1)=((1,1),D_k^2,D_{k-2}^1,D_{k-4}^2,D_{k-6}^1,\ldots,D_3^2,(2,2),\ldots,(k+1,k+1),$$
$$D_{k-1}^1,D_{k-3}^2,\ldots, D_4^2,D_2^1,(1,\ell),(n-k+2,n-k+2),\dots,(\ell+k-1,\ell+k-1),$$
$$D_{k-1}^2,D_{k-3}^1,\ldots,D_{4}^1,D_2^2,D_{k}^1,D_{k-2}^2,\ldots,D_{5}^2,D_3^1,(\ell,\ell),\ldots,(k,k)).$$
It is clear that $D_i\subseteq L(1,1)$ for any $i\neq1$. Suppose by way of contradiction that there exists an element $(j,j)\in D_1$, with $(j,j)\not \in L(1,1)$.
Hence $L(j,j)\cap L(1,1) =\emptyset$, which implies $L(j,j)\subseteq D_1$ which is absurd since $\gcd(n,k-1)=2$ and $2,\ell\in E$.
\endproof

\begin{ex}\label{ExampleAlmost}
Consider a cyclically almost $7$-diagonal array of size $14$ in standard form with the extra filled position $(1,5)$.
We show in the table below the tour obtained in the proof of Proposition \ref{3mod4}, where cells are tagged as in previous examples.

$$\begin{array}{r|r|r|r|r|r|r|r|r|r|r|r|r|r|r|r|}     &\downarrow & \uparrow   & \downarrow  & \downarrow& \uparrow& \downarrow   &\downarrow  &\downarrow   &\downarrow    &\downarrow    &\downarrow    &\downarrow   & \downarrow  &\downarrow \\\hline
 \rightarrow    &0 &    &   & &50 &    &  &   &   5  & 59   &  88  & 43  & 23  & 77 \\\hline
 \rightarrow  & 78 & 26 &    &   & & &   &  & &    6  & 60   &  89  & 44  &  24 \\\hline
 \rightarrow	&25 &	 47 & 52 &    & &   &  &  & & &    7  & 61   &  90  & 45  \\\hline
 \rightarrow&46 &	93 &	 48 & 28 &    & &   &  &  & &  & 8  & 62   &  91     \\\hline
 \rightarrow&92  & 65   &  94  & 49& 96 &  &   &  &  &   & & &     9 & 63  \\\hline
 \rightarrow&64 & 12   &  66  & 95 &	 68 & 30 &    &&   &  &   & &  & 10   \\\hline
 \rightarrow&11  & 33   &  13  & 67 &	15 &	 69 & 98 & &   &  &   &   & &   \\\hline
 \rightarrow&&79  & 34  &  14  & 36 &	16 &	 70 & 32  &   &  &   &   & &   \\\hline
 \rightarrow&&&80  & 35   &  82  & 37 &	17 &	 71 & 51   &  &   &   & &   \\\hline
 \rightarrow&  & &  &81  & 54   & 83  & 38 &	18 &	72 & 27 &    &   & &   \\\hline
 \rightarrow& & &  &  &1  & 55   &  84  & 39 &	19 &	 73 & 53    &   & &   \\\hline
 \rightarrow& &   &&  &   &2  & 56   &  85  & 40 &	20 &	74 & 29  & &   \\\hline
 \rightarrow& &   &&  &   &   &3  & 57   &  86  & 41 &	21 &	75 & 97  &   \\\hline
 \rightarrow& &   &  &&   &   & & 4  & 58   &  87  & 42 &	22 &	 76 & 31   \\\hline
			\end{array}$$
\end{ex}
Now we investigate the case $k$ even. We recall that, in this case, in order to have a solution the size $n$ of an almost $k$-diagonal
array can be both even and odd.
\begin{prop}
 Let $A$ be a cyclically almost $2$-diagonal
array of size $n\geq3$ with extra filled cell $(1,\ell)$.
The vectors $R:=(1,\ldots,1)$ and $C:=(1,\ldots,1,c_\ell=-1,1,\ldots,1)$ are a solution of $P(A)$.
\end{prop}
\proof
It is sufficient to note that $L(1,1)=((1,1),(\ell+1,\ell),(\ell+2,\ell+1),(\ell+3,\ell+2),\ldots,(n,n-1),(1,n),(2,1),
(3,2),\ldots,(\ell,\ell-1),(1,\ell),(n,n),(n-1,n-1),\ldots,(2,2))$.
\endproof

For $k\geq4$, we start by deriving some properties that must be satisfied by  any solution.
\begin{prop}\label{necckeven}
  Let $k\geq4$ be an even integer and let
$A$ be a cyclically almost $k$-diagonal
array of size $n> k$ whose extra filled position is $(1,\ell)$.
If the vectors $R:=(1,\ldots,1)$ and $C\in \{-1,1\}^n$ are a solution of $P(A)$ then
\begin{itemize}
  \item[1)] $c_\ell=-1$;
  \item[2)] $\gcd(|E|,\frac{k}{2})=1$;
  \item[3)] denoted by $d=\gcd(n,k-1)$, the list $E$ covers all the congruence classes modulo $d$ if $\ell\equiv1\pmod d$,
otherwise the list $E$ covers all the congruence classes modulo $d$ except at most the class of $1$.
\end{itemize}
\end{prop}
\proof
1) If, by way of contradiction, we suppose that $c_\ell=1$, then $L(1,1)\subseteq D_1\cup D_3 \cup \ldots \cup D_{k-1}\cup \{(1,\ell)\}$, hence
$L(1,1)$ does not cover $F(A)$.

2) Note that since $k\geq4$, $(\ell+2,\ell-1)\in D_4$ is a filled position
hence we can consider the list $L(\ell+2,\ell-1)$.
If $(1,\ell)\not\in L(\ell+2,\ell-1)$ then $R$ and $C$ are not a solution.
So, suppose that $(1,\ell)\in L(\ell+2,\ell-1)$.
Since, by 1), $c_\ell=-1$ we have that $S_{R,C}((i,j))=(1,\ell)$ if and only if $(i,j)=(\ell,\ell-1)$ .
Also, note that applying $S_{R,C}$ to a filled cell of the diagonal $D_{2a}$ and of the $j$-th column we obtain a cell of the diagonal $D_{2a-2}$ (modulo $k$)
and of the $(j+1)$-th column (modulo $n$) or a cell of $D_{2a}$.
Hence we have $S_{R,C}^\alpha((\ell+2,\ell-1))=(\ell,\ell-1)\in D_2$ implies that $\alpha=\lambda n$ for a suitable $\lambda$.
We also note that given $(i,j)\neq(\ell,\ell-1)$ and $(i,j)\in D_h$ with $h=2,\ldots,k$, $S_{R,C}((i,j))\in D_h$ if $j+1\not\in E$,
otherwise $S_{R,C}((i,j))\in D_{h-2}$ where the subscripts are taken modulo $k$.
Hence, after having applied $\alpha=\lambda n$ times $S_{R,C}$ to $(\ell+2,\ell-1)$  we are in the diagonal of index $4-2\lambda|E|\pmod k$,
but we have already observed that $S_{R,C}^\alpha((\ell+2,\ell-1))=(\ell,\ell-1)\in D_2$. So we obtain $4-2\lambda|E|\equiv2\pmod k$, that is
$\lambda|E|\equiv 1\pmod{ \frac{k}{2}}$, which implies that $\gcd(|E|,\frac{k}{2})=1$.

3) By way of contradiction, we suppose that there exists a congruence class $[i]$ modulo $d$, with $i\not\equiv1 \pmod d$
which is not covered by $E$.
Let $(f,f)\in \mathcal{C}_i$ and note that  $S_{R,C}((f,f))\in \mathcal{C}_i$ since $(f,f)\neq(1,1)$
because $i\not\equiv1\pmod d$ and since $E$ does not cover the congruence class of $i$.
Hence $L(f,f)\subseteq \mathcal{C}_i$. So $L(1,1)$ does not cover all the filled position of $A$.
\endproof

\begin{prop}
  Let $k\geq 4$ be an even integer and let
$A$ be a cyclically almost $k$-diagonal
array of size $n> k$ whose the extra filled position is $(1,\ell)$. If $\gcd(n,k-1)=1$ then
the vectors $R:=(1,\ldots,1)$ and $C:=(1,\ldots,1,c_\ell=-1,1,\ldots,1)$ are a solution of $P(A)$.
\end{prop}
\proof
Let $R$ and $C$ be as in the statement.
For any $i=2,\ldots,k$ we set
$D_i^1=\{(x,y)\in D_i\ |\ 1\leq y\leq \ell-1\}$ and $D_i^2=D_i\setminus D_i^1$.
It is easy to see that
$$L(1,1)=((1,1),D_k^2,D_{k}^1,D_{k-2}^2,D_{k-2}^1,\ldots, D_2^2,D_2^1,(1,\ell),(n-k+2,n-k+2),\ldots,(k,k)).$$
Set $B=((n-k+2,n-k+2),\ldots,(k,k))\subseteq L(1,1)$.
Note that $(n-k+2,n-k+2)\in D_1$ hence either $(\ell+k-1,\ell+k-1)\in B$ or $B\subseteq D_1$.
It is not hard to see that $B\subseteq D_1$ implies that
$B=(b_1,b_2,\ldots,b_{|B|})$ where
$b_\alpha=(n-k+2-(\alpha-1)(k-1),n-k+2-(\alpha-1)(k-1))$ for $\alpha=1,\ldots,|B|$.
But $n-k+2-(\alpha-1)(k-1)\equiv k\pmod n$ if and only if $\alpha\equiv n-1\pmod n$, which implies $|B|=n-1$ that is $B=D_1\setminus{(1,1)}$. Since $2\leq \ell\leq n-k+1$ we have $\ell+k-1\not\equiv 1 \pmod n$.
So we can suppose that $(\ell+k-1,\ell+k-1)\in B$ which implies that
$$L(1,1)=((1,1),D_k^2,D_{k}^1,D_{k-2}^2,\ldots,D_2^2,D_2^1,(1,\ell),(n-k+2,n-k+2),\ldots,
$$
$$(\ell+k-1,\ell+k-1),D_{k-1}^2,D_{k-1}^1,D_{k-3}^2,\ldots,D_3^2,D_3^1,(\ell,\ell),\dots,(k,k)).$$
It is clear that $D_i\subseteq L(1,1)$ for any $i\neq1$. Suppose, by way of contradiction, that there exists an element $(j,j)\in D_1$, with $(j,j)\not \in L(1,1)$.
Hence $L(j,j)\cap L(1,1) =\emptyset$ implying $L(j,j)\subseteq D_1$, which is absurd since $\gcd(n,k-1)=1$.
\endproof

\section{A recursive construction}\label{SecRecursive}
We have already seen in Theorem \ref{prop:necc} that if $A$ admits nontrivial closed subarrays then $P(A)$ has no solution.
Note that if $A$ has exactly $t$ minimal closed subarrays, adding to $A$ less than $t-1$ extra filled cells, we get again an array with nontrivial closed subarrays.
Hence in order to obtain a solution of the problem for ``arrays with almost closed subarrays'' it is necessary to add at least $t-1$  filled cells.
We start investigating the case in which $t=2$ and the arrays are placed diagonally.

\begin{prop}\label{blockdiag}
  Let $B_1$ and $B_2$ be two arrays with no empty row and no empty column.
  Let $A$ be an array of the form  $$A=\begin{array}{|r|r|}\hline
    B_1 & E_1 \\ \hline
    E_2 & B_2\\ \hline
  \end{array}$$
  where $E_1\cup E_2$ has exactly one filled position. Then $P(A)$ has a solution if and only if both $P(B_1)$ and $P(B_2)$ do.
\end{prop}
\proof
Let $(i,j)$ be the unique filled position of $E_1\cup E_2$.
Firstly, we suppose that $(i,j)\in E_2$.

Let $R_t$ and $C_t$ be a solution of $P(B_t)$ for $t=1,2$.
We will show that $R:=(R_1,R_2)$ and $C:=(C_1,C_2)$ are a solution of $P(A)$.
Note that the element $(i',j')$ such that $S_{R,C}((i',j'))=(i,j)$ belongs to $B_1$ and that $S_{R,C}((i,j))\in B_2$.
Hence the list $L(i,j)$ has elements both in $B_1$ and in $B_2$.
Suppose now that there exists $(f,g)\in F(A)\setminus L(i,j)$ and that $L(f,g)$ contains elements both of $B_1$ and $B_2$. Since for any $(i_1,j_1)\in B_1$, $S_{R,C}((i_1,j_1))\not\in B_2$, $(i,j)\in L(f,g)$ which is absurd.
Hence, $L(f,g)$ contains only elements of $B_t$ for some $t\in\{1,2\}$.
Since $R$ and $C$ restricted to $B_t$ are nothing but $R_t$ and $C_t$, which solve $P(B_t)$,
the list $L(f,g)$ has to contain all the elements of $B_t$. We deduce that
 $L(f,g)\cap L(i,j)\neq \emptyset$, which is absurd.

Suppose now that $R$ and $C$ are a solution of $P(A)$, namely that $L(i,j)$ covers all the filled positions of $A$.
Since, as already remark, $S_{R,C}((i,j))\in B_2$ and the element $(i',j')$ such that $S_{R,C}((i',j'))=(i,j)$ belongs to $B_1$,
we have that $L(i,j)=((i,j), \bar{B_2},\bar{B_1})$ where $\bar{B_t}$ is a list of filled cells of $B_t$, for $t=1,2$.
The list $L(i,j)$ covers all the filled positions of $A$; hence $\bar{B_t}=F(B_t)$ for $t=1,2$, so the restrictions of $R$ and $C$ to $B_t$ give a solution of $P(B_t)$.

 If $(i,j)\in E_1$ the proof can be done in a similar way.
\endproof

\begin{ex}\label{exblockdiag}
Let $B_1$ be a cyclically $3$-diagonal array of size $7$ in standard form and let $B_2$ be a totally filled $3\times 4$ array.
A solution of $P(B_1)$ and one of $P(B_2)$ are presented in Example  \ref{ex5} and in Theorem \ref{rettangolopieno}, respectively.
Hence, joining these solutions, we have that $R:=(-1, 1, 1, 1, 1, 1, 1 ,1,1,1)$ and
$C:=(-1, -1, 1, 1, 1, -1, -1,1,1,1,-1 )$ are a solution of $P(A)$, for the array $A$ shown in the table below, where cells are tagged as usual.
\begin{center}
$A:=\begin{array}{r|r|r|r|r|r|r|r||r|r|r|r|}
   &\uparrow & \uparrow   & \downarrow   &  \downarrow  & \downarrow   &\uparrow   & \uparrow & \downarrow   & \downarrow   &  \downarrow  & \uparrow \\\hline
\leftarrow &31 &    &    &    &    &  25 &  28&    &    &    &  \\\hline
\rightarrow &14 &  27 &       &    &   &  &  30 &    &    &    & \\\hline
\rightarrow  &26 & 21 &  13&     &    &  &  &    &    &    & \\\hline
\rightarrow  &&  15 &  22 &  20 &       &    & &    &    &    &   \\\hline
\rightarrow    & &    &   16 &   23 &   33 &    & &    &    &    &    \\\hline
\rightarrow     &&     &   &  17 & 24 & 19 &    &    &    &    & \\\hline
\rightarrow   &    &  &   &  &  18 &  29 &  32 &    &    &    & \\\hline
\hline
\rightarrow  &    & &  &   &    &     &  &  9   & 2   &   7  &  4   \\\hline
\rightarrow  &    &  &0 &   &    &     & &  5   &  10   &   3  &  12   \\\hline
\rightarrow  &    &   &&   &    &     & &  1   & 6   &  11  &  8   \\\hline
  \end{array}$
\end{center}
\end{ex}
We have to point out that in Proposition \ref{blockdiag} it is not necessary to require that the arrays $B_1$ and $B_2$ are placed diagonally.
In fact the same proof holds also when the rows and the columns of $B_1$ (resp. $B_2$) are not consecutive in $A$,
as shown in the following example.

\begin{ex}
Let $B_1=\mathcal{R}_1 \cap \mathcal{C}_1$ where $\mathcal{R}_1=(R_1,R_2,R_3,R_4,R_5,R_7,R_8)$
and $\mathcal{C}_1=(C_1,C_2,C_3,C_4,\\C_5,C_8,C_9)$ and set
$B_2=\mathcal{R}_2 \cap \mathcal{C}_2$ where $\mathcal{R}_2=(R_6,R_9,R_{10})$
and $\mathcal{C}_2=(C_6,C_7,C_{10},C_{11})$.
Note that $B_1$ and $B_2$ are the same arrays of Example \ref{exblockdiag}, but now their rows and columns are not consecutive in $A$.
Also in this case starting from the solutions of $P(B_1)$ and $P(B_2)$ it is possible to obtain a solution of $P(A)$ as shown below:
\begin{center}
$A:= \begin{array}{r|r|r|r|r|r|r|r|r|r|r|r|}
   &\uparrow & \uparrow   & \downarrow   &  \downarrow  & \downarrow    & \downarrow   & \downarrow  &\uparrow   & \uparrow &  \downarrow  & \uparrow \\\hline
\leftarrow &31 &    &    &    &    &    &    &25 &  28&    &  \\\hline
\rightarrow &14 &  27 &       &    &    &    &&  &  30    &    & \\\hline
\rightarrow  &26 & 21 &  13&     &    &  &  &    &    &    & \\\hline
\rightarrow  &&  15 &  22 &  20 &       &    & &    &    &    &   \\\hline
\rightarrow    & &    &   16 &   23 &   33 &    & &    &    &    &    \\\hline
\rightarrow  &    & &  &   &    &    9   & 2  &  &    &   7  &  4   \\\hline
\rightarrow     &&     &   &  17 & 24   &    && 19 &      &    & \\\hline
\rightarrow   &    &  &   &  &  18 &    &  &  29 &  32   &    & \\\hline
\rightarrow  &    &  &0 &   &    &      5   &  10& &    &   3  &  12   \\\hline
\rightarrow  &    &   &&   &    &      1   & 6& &    &  11  &  8   \\\hline
  \end{array}$
\end{center}
\end{ex}

\begin{rem}
Let $A$ be an array with $t$ minimal closed subarrays $B_i$, for $i=1,\ldots,t$, such that $P(B_t)$ has a solution.
Iterating the previous procedure, it is possible to add $t-1$ extra filled positions (in a suitable way) to $A$ in order to obtain a new array $A'$ such that $P(A')$ has a solution.
\end{rem}

\section*{Aknowledgments}
This research was partially supported by Italian Ministry of Education, Universities and Research under Grant PRIN 2015 D72F16000790001 and by INdAM-GNSAGA.

\end{document}